\documentclass[10pt]{article}
\usepackage{url}
\usepackage{amsthm,amssymb}
\usepackage{amsfonts}
\usepackage{verbatim}
\usepackage{setspace}
\usepackage{float}
\theoremstyle{plain}
\usepackage[ruled,noend]{algorithm2e}
\newtheorem{theorem}{Theorem}

\newtheorem{proposition}[theorem]{Proposition}

\theoremstyle{definition}

\DeclareSymbolFont{AMSb}{U}{msb}{m}{n}
\DeclareMathSymbol{\N}{\mathbin}{AMSb}{"4E}
\DeclareMathSymbol{\Z}{\mathbin}{AMSb}{"5A}
\DeclareMathSymbol{\R}{\mathbin}{AMSb}{"52}
\DeclareMathSymbol{\Q}{\mathbin}{AMSb}{"51}
\DeclareMathSymbol{\I}{\mathbin}{AMSb}{"49}
\DeclareMathSymbol{\C}{\mathbin}{AMSb}{"43}

\title{Nontorsion Points of Low Height on Elliptic Curves over Quadratic Fields}
\author{Graeme Taylor\footnote{Heilbronn Institute for Mathematical Research, University of Bristol, University Walk, Bristol BS8 1TW, UK}}
\date{}
\begin{document}
\maketitle
\begin{abstract}\noindent We give examples of points with particularly low height on elliptic curves over quadratic fields, recovered by a search over elliptic divisibility sequences. The smallest example identified satisfies $d\hat{h}(P)=0.0077127\ldots$: improving on the previous smallest for curves over quadratic fields of $d\hat{h}(P)=0.0194426\ldots$ given in \cite{EverestWard}; and comparable with some of the examples of smallest height on curves over $\Q$ tabulated in \cite{Elkies}. \end{abstract}

\section{Introduction}
The \emph{elliptic Lehmer problem} is concerned with lower bounds on the height of points on elliptic curves; more precisely, if $E$ is an elliptic curve over $K$ a number field of degree $d$, then it is conjectured that for any nontorsion point $P\in E(K)$ the canonical height $\hat{h}(P)$ is bounded below by $c/d$ for $c>0$ a uniform constant. Elkies maintains a website \cite{Elkies} which includes a table of the smallest known examples over $\Q$: there are $54$ examples with $\hat{h}<0.01$; the best being the point $(0,0)$ of height $0.0044571\ldots$ on $y^2+253xy+119700y=x^3+5320x^2$.\\ 

No analogous tables for general number fields exist, but \cite{EverestWard} contains the examples $(3-12\omega,-108\omega^2)$ of height $0.01032\ldots$ on $y^2=x^3-243x+10368\omega+3726$, where $\omega$ is a non-trivial cube root of unity; and $(3-9u,108-108u)$ of height $0.00971\ldots$ on $y^2=x^3+(1215u-2214)x+40878-23328u$, where $u=(1+\sqrt{5})/2$. The elliptic Lehmer problem implies a lower bound on the growth rate of elliptic divisibility sequences, and Everest's and Ward's examples were recovered from such sequences with arithmetically simple initial terms, which they suggest as a promising search space.\\ 

This note gives the results of such a search for points of small height, for elliptic curves over some quadratic fields. Two examples with height less than $0.005$ (and thus competitive with the entries in Elkies' tables) were found: the point $(497\sqrt{3}+860,7523\sqrt{3}+13029)$ of height $0.0038563\ldots$ on $y^2+\sqrt{3}xy + y =x^3-\sqrt{3}x^2 -(1544879\sqrt{3}+2675806)x +1407381720\sqrt{3} +2437656645$; and $(2\omega-4,3\omega+6)$ of height $0.0047223\ldots$ on $y^2+xy=x^3-2\omega x^2 -(\omega+2)x +5\omega+2$, where $\omega=\frac{1}{2}+\frac{\sqrt{-7}}{2}$. In all, $30$ examples with $\hat{h}<0.01$ are given in Table \ref{lowheighttable} in Section 4; these are preceeded by a summary of the relevant theory in Section 2, and a description of techniques for efficient implementation of the search in Section 3.

\section{Elliptic Divisibility Sequences and Elliptic Curves}

An \emph{elliptic divisibility sequence} is a sequence of integers ${u_n}:u_0,u_1,u_2,\ldots,u_i,\ldots$ which is a particular solution of \begin{equation}\label{edsrecursion}w_{m+n}w_{m-n}=w_{m+1}w_{m-1}w_n^2-w_{n+1}w_{n-1}w_m^2\end{equation} such that $u_n$ divides $u_m$ whenever $n$ divides $m$. 
A \emph{proper} elliptic divisibility sequence additionally satisfies \begin{equation}\label{properEDS} u_0=0, u_1=1, u_2u_3\neq0.\end{equation} By \cite{Ward}, these are the only elliptic divisibility sequences of arithmetic interest, and a tuple $u_2,u_3,u_4$ defines a proper elliptic divisibility sequence if and only if they are integral with $u_2|u_4$.

Let $K$ be a number field of degree $d$ with ring of integers $\mathcal{O}_K$. For $a_1,a_2,a_3,a_4,a_6\in\mathcal{O}_K$ an elliptic curve $E$ over $K$ is given by the generalized Weierstrass equation \begin{equation} \label{generalizedWeierstrass}y^2+a_1xy+a_3y=x^3+a_2x^2+a_4x+a_6.\end{equation}  Defining quantities $b_2,b_4,b_6,b_8$ as 
\begin{eqnarray*}b_2&=&a_1^2+4a_2,\\
b_4&=&2a_4+a_1a_3,\\
b_6&=&a_3^2+4a_6,\\
b_8&=&a_1^2a_6+4a_2a_6-a_1a_3a_4+a_2a_3^2-a_4^2,\end{eqnarray*}
we obtain the \emph{discriminant} $\Delta$ of $E$, 
\begin{equation}\label{discriminant} \Delta=-b_2^2b_8-8b_4^3-27b_6^2+9b_2b_4b_6\in\mathcal{O}_K,\end{equation}
and may construct polynomials $\psi_n\in\mathcal{O}_K[x,y]$ by
\begin{eqnarray*}
\psi_0&=&0,\\
\psi_1&=&1,\\
\psi_2&=&2y+a_1x+a_3,\\
\psi_3&=&3x^4+b_2x^3+3b_4x^2+3b_6x+b_8,\\
\psi_4&=&\psi_2(2x^6+b_2x^5+5b_4x^4+10b_6x^3+10b_8x^2+(b_2b_8-b_4b_6)x+b_4b_8-b_6^2)\\
\mbox{and then inductively for $n\ge3$ by}&&\\ 
\psi_{2n-1}&=&\psi_{n+1}\psi_{n-1}^3-\psi_{n-2}\psi_n^3,\\
\psi_{2n}&=&\psi_n(\psi_{n+2}\psi_{n-1}^2-\psi_{n-2}\psi_{n+1}^2)/\psi_2.
\end{eqnarray*}
The zeros of $\psi_n$ are the $x$-coordinates of the points on $E$ with order dividing $n$. Let $\psi_n(P)$ denote $\psi_n$ evaluated at the point $P=(x,y)$. Then, for any $K$-rational point $P=(x,y)$ on $E$, the sequence ${u_n=\psi_n(P)}$ satisfies (\ref{edsrecursion}) and (\ref{properEDS}); that is, for any $P\in E(K)$ we have a corresponding proper elliptic divisibility sequence. If $P$ is not a torsion point, then the terms of this sequence are always non-zero.\\

Let $M_K$ be the set of valuations of $K$, each corresponding to some absolute value $|\cdot|_v$, with $K_v$ the corresponding completion of $K$. The \emph{na\"{\i}ve height} $h(\alpha)$ of $\alpha\in K$ is
\[h(\alpha)=\frac{1}{d}\sum_{v\in M_K} \mbox{log max}\{1,|\alpha|_v\},\] and for a finite point $P\in E(K)$ we define $h(P)=h(x(P))$, with $h(P_0)=0$ for $P_0$ the point at infinity. Then the \emph{global canonical height} is the function $\hat{h}:E(K)\rightarrow \R$ given by
\begin{equation}\label{tateheight}\hat{h}(P)=\frac{1}{2}\lim_{n\rightarrow\infty}4^{-n}h([2^n]P).\end{equation}
 
(We note that there are two competing definitions of $\hat{h}(P)$ in the literature, differing by a factor of $2$; in (\ref{tateheight}) and throughout we follow \cite{Silverman} in taking the smaller, which is consistent with \cite{Elkies}. The alternative is more natural in the context of the Birch and Swinnerton-Dyer conjecture and, in particular, is the value returned by the height functions in computer algebra systems such as Magma \cite{magma} or SAGE \cite{sage}.)

Let $D=N_{K|\Q}(\Delta)$ and $T$ the set of rational primes which divide $D$. Then, given an algebraic integral point $P$ on $E(K)$ we have (by Theorem 1 of \cite{EverestWard})
\begin{equation}\label{everestwardheight} \hat{h}(P)=\displaystyle\frac{1}{d}\lim_{n\rightarrow\infty}\frac{1}{n^2}\log{F_n},\end{equation}
where \[F_n=E_n\prod_{p\in T} |E_n|_p  \mbox{ for }  E_n=|N_{K|\Q}(\psi_n(P))|.\]

Further, it is usually sufficient in practice to compute (\ref{everestwardheight}) via the gcd of $E_n$ and $E_{n+1}$, as in formula $(21)$ of \cite{EverestWard}:
\begin{equation}\label{everestwardheightbygcd} \hat{h}(P)=\displaystyle\frac{1}{d}\lim_{n\rightarrow\infty}\frac{1}{n^2}\log{\left(\frac{E_n}{\mbox{gcd}(E_n,E_{n+1})}\right)} \approx \tilde{h}_n=\frac{1}{dn^2}\log{\left(\frac{E_n}{\mbox{gcd}(E_n,E_{n+1})}\right)} \end{equation}
and, to a few significant figures, good approximations $\tilde{h}_n$ are obtained for values of $n$ as small as $128$. In the following section, we will describe efficient procedures for recovering pairs $E_n,E_{n+1}$ and hence $\tilde{h}_n$.

\begin{proposition}\label{conjugatefactors} The proper elliptic divisibility sequences defined by $u_2,u_3,u_4$ and $\overline{u_2},\overline{u_3},\overline{u_4}$ correspond to points $P$ and $P'$ of the same height.\end{proposition} 

We may also perform the reverse procedure, of constructing a point on an elliptic curve from an elliptic divisibility sequence. The following unpleasant formulae from \cite{Ward} give a point $(x,y)$ on an elliptic curve of form $y^2=4x^3-g_2x-g_3$ from the proper elliptic divisibility sequence defined by the tuple $u_2,u_3,u_4$:
\[ x= \displaystyle \frac{u_4^2+2u_2^5u_4+4u_2u_3^3+u_2^{10}}{12u_2^4u_3^2}, y=-u_2\]
\[g_2=\frac{u_2^{20}+4u_2^{15}u_4-16u_2^{12}u_3^3+6u_2^{10}u_4^2-8u_2^7u_3^3u_4+4u_2^5u_4^3+16u_2^4u_3^6+8u_2^2u_3^3u_4^2+u_4^4}{12u_2^8u_3^4}\]
\begin{eqnarray*}g_3&=&\frac{-1}{216u_2^{12}u_3^6}(u_2^{30}+6u_2^{25}u_4-24u_2^{22}u_3^3+15u_2^{20}u_4^2-60u_2^{17}u_3^3u_4\\
&&+20u_2^15u_4^3+120u_2^{14}u_3^6-36u_2^{12}u_3^3u_4^2+15u_2^{10}u_4^4-48u_2^9u_3^6u_4\\
&&+12u_2^7u_3^3u_4^3+64u_2^6u_3^9+6u_2^5u_4^5+48u_2^4u_3^6u_4^2+12u_2^2u_3^3u_4^4+u_4^6)\end{eqnarray*}

\begin{proposition}\label{firstfactorpositive} The proper elliptic divisibility sequences defined by $u_2,u_3,u_4$ and $-u_2,u_3,-u_4$ give points $P=(x,y)$ and $(x,-y)=-P$ on the same curve.\end{proposition} 

In practice, we may often use the much simpler formulae given in \cite{Shipsey} to recover the point $(0,0)$ on a curve of form (\ref{generalizedWeierstrass}), where the $a_i$ are given by either
\begin{equation}\label{shipsey1} a_1=0, a_3=u_2, a_4=\frac{u_4+u_2^5}{2u_2u_3}, a_2=\frac{u_3+a_4^2}{u_2^2}, a_6=0,\end{equation} or
\begin{equation}\label{shipsey2} a_1=1, a_3=u_2, a_4=\frac{u_4-u_2^2u_3+u_2^5}{2u_2u_3}, a_2=\frac{u_3+a_1a_3a_4+a_4^2}{u_2^2}, a_6=0.\end{equation}

We may therefore search for points of low height by: generating candidate tuples $u_2,u_3,u_4$ that satisfy (\ref{edsrecursion}) and (\ref{properEDS}) (for which it suffices to take the $u_i$ integral with $u_2|u_4$); obtaining an estimate $\tilde{h}_n$; and, when $\tilde{h}_n$ is sufficiently small, recovering a corresponding point and curve $P,E$ with $\hat{h}(P)\approx\tilde{h}_n$.

\section{Efficient Computation with Elliptic Divisibility Sequences}
We wish to rapidly compute the terms $u_n,u_{n+1}$ of a proper elliptic divisibility sequence. The inductive construction of polynomials $\psi_n$ for $n$ around $2k$ from polynomials with index around $k$ suggests a double-and-add addition chain approach; in \cite{Shipsey} this is made precise, and described in full generality. For our application, a slightly simpler description can be given, since instead of arbitrary chains, we can restrict to a sequence of doublings.

We define the septuple $\langle u_k\rangle$ to be the terms $\{u_{k-3},u_{k-2},u_{k-1},u_k,u_{k+1},u_{k+2},u_{k+3}\}$. Given $\langle u_k\rangle$, we can then recover $\langle u_{2k}\rangle$, since from (\ref{edsrecursion}) we have 
\begin{eqnarray} u_{2l}&=&u_l(u_{l+2}u_{l-1}^2-u_{l-2}u_{l+1}^2)/u_2\label{evendouble}\\
u_{2l+1}&=&u_{l+2}u_l^3-u_{l-1}u_{l+1}^3\label{odddouble} \end{eqnarray}
and setting $l=k-1,k,k+1$ in (\ref{evendouble}), $l=k-2,k-1,k,k+1$ in (\ref{odddouble}) generates the desired terms from elements of $\langle u_k\rangle$ and $u_2$ only. 

For a proper elliptic divisibility sequence we know $u_0,\ldots,u_4$; from (\ref{edsrecursion}) we have $u_{-n}=-u_n$ for any $n\in\N$, so $u_{-1}=-1$; and from (\ref{odddouble}) we have $u_5=u_4u_2^3-u_3^3$. So we may construct $\langle u_2\rangle$ and hence, iteratively, $\langle u_n\rangle$ for any $n=2^m$, from which we may recover  $u_n,u_{n+1}$ to compute $\tilde{h}_n$.\\

In \cite{Stange}, elliptic divisibility sequences are generalised to elliptic nets, and efficient algorithms are given for the double/double-and-add operation on a 'block' centred at $k$, which contains $\langle u_k\rangle$ as a subset. Specialising these to elliptic divisibility sequences, we note that a speed-up in the function $double: \langle u_k\rangle \mapsto \langle u_{2k}\rangle$ can be gained by pre-computing the frequently used expressions $A_i=u_i^2$ and $B_i=u_{i-1}u_{i+1}$ for $i=k-2,\ldots k+2$. Assuming $u_2^{-1}$ is pre-computed for a given sequence, we can then perform $double$ as described in Algorithm \ref{doublealgorithm}, which reduces the operation count from 14 squarings and 28 multiplications to 5 squarings and 22 multiplications.

\begin{algorithm}[h]
\KwIn{$[V_1,V_2,V_3,V_4,V_5,V_6,V_7]=\langle u_k \rangle, u_2^{-1}$}
\KwOut{$\langle u_{2k} \rangle$}
\For{$i=1,\ldots,5$}
{
	$A[i]=V[i+1]^2$\;
	$B[i]=V[i]V[i+2]$\;	
}
\For{$i=0\ldots,3$}
{
	$V[2i+1]=B[i+2]A[i+1]-B[i+1]A[i+2]$\;
	\uIf{$i>0$}{$V[2i]=(B[i+2]A[i]-B[i]A[i+2])u_2^{-1}$}
}		
\Return{$V$}
\caption{$double$}
\label{doublealgorithm}
\end{algorithm}

To exhaustively search over a selection of tuples we then proceed as in Algorithm \ref{searchalgorithm}. Note that we ensure $u_2|u_4$; and pre-compute $u_2^{-1}$ once per $u_2$ value, rather than once per sequence (which reduces runtime by about a third in practice). Tuples with $E_n$ or $E_{n+1}$ equal to zero necessarily correspond to torsion points, and thus can be discarded. Any tuple amongst those returned by Algorithm \ref{searchalgorithm} is then passed to secondary testing: first we attempt to construct the corresponding curve, and if this is nonsingular, we check the order of the corresponding point is infinite, then confirm its height with both an estimate $\tilde{h}_n$ for a larger $n$, and a direct computation of $\hat{h}(P)$ in Magma\footnote{We note that some of the recovered points cause difficulties for the height calculation in SAGE, even with precision set as high as 2048-bit; but in each case Magma confirms the estimates given by (\ref{everestwardheightbygcd}).}. 

\begin{algorithm}[h]
\KwIn{Search spaces $S_1,S_2,S_3$; Height bound $H$, Iteration bound $I$}
\KwOut{Defining terms for eds with height estimate $0<\tilde{h}<H$}
$Secondary=List([])$\; 
\For{$\alpha\in S_2$}
	{
		$u_2=\alpha$\;
		$iu_2=\alpha^{-1}$\;
		\For{$\beta\in S_2$}
			{
			$u_3=\beta$\;
			\For{$\gamma\in S_3$}
				{
				$u_4=\alpha\gamma$\;
					$V=[-1,0,1,u_2,u_3,u_4,u_4u_2^3-u_3^3]$\;
					$n=2$\;
					\For{$i=1,\ldots, I$}
						{
						$V=double(V,iu_2)$\;
						$n=2n$\;	
						}
					$E_n=|N_{K|\Q}(V_4)|$\;
					$E_{n+1}=|N_{K|\Q}(V_5)|$\;
					\uIf{$E_n\neq0$  and $E_{n+1}\neq0$ }
					  {
					  $\tilde{h}=\log(E_n/gcd(E_n,E_{n+1}))/dn^2$\;
					  \uIf{$0<\tilde{h}<H$}{$Secondary.append([u_2,u_3,u_4])$}
					  }
				}
			}
	}
\Return{$Secondary$}
\caption{$search$}
\label{searchalgorithm}
\end{algorithm}

\section{Results}

We consider quadratic fields $K=\Q(\sqrt{D})$ for $D\in\{-7,-6,-5,-3,-2,-1,2,3,5,6,7\}$. For such fields we have $\mathcal{O}_K=\Z[\omega]$ where $\omega=\frac{1}{2}+\frac{\sqrt{D}}{2}$ if $D\equiv 1$ modulo $4$, or $\omega=\sqrt{D}$ for $D\equiv 2$ or $3$. We seek to test the tuples of form $u_2=\alpha$, $u_3=\beta$, $u_4=\alpha\gamma$ for \[\alpha,\beta,\gamma\in S:=\left\{x+y\omega\,|\,x,y\in\Z, |x|+|y|>0, \max{\{|x|,|y|\}}\le c\right\},\] for some bound $c$. By Proposition \ref{firstfactorpositive} we may restrict to $\alpha\in S^+$, where 
\[ S^+=\left\{x+y\omega\,|\,x\in\{1,\ldots,c\}, y\in\{-c\ldots c\}\right\}\cup \{y\omega\,|\,y\in\{1,\ldots,c\}\}\] For $c=9$ this gives a total of $23,328,000$ tuples to test for each $D$; with $H=0.01$ and $I=6$ (so we estimate $\hat{h}$ by $\tilde{h}_{128}$), runtimes for Algorithm \ref{searchalgorithm} in SAGE were in the range $12.5$-$20$cpuhours\footnote{where a `CPU' is one core of an Intel X5650 at 2.67ghz; parallelisation across cores is trivial by partitioning $S^+$.} per $D$. Based on timings for $c=3$ and $D=-7$ (the choice of $D$ for which Algorithm \ref{searchalgorithm} was slowest), this represents a speedup of around a factor of $50$ over generating $P$ for each tuple and using SAGE's built-in height function to compute $\hat{h}(P)$, despite being \emph{more} accurate\footnote{although we note that performance degrades substantially for $c\ge 10$, which is why we have restricted to $c=9$.}. 

Representatives of the points $P$ on curves $E$ satisfying $\hat{h}(P)<0.01$ found by Algorithm \ref{searchalgorithm} with $S_1=S^+,S_2=S_3=S$, and the tuple $h_2,h_3,h_4$ by which they were identified, are given in Table \ref{lowheighttable}. In the interests of brevity, if points $P,E$ and $P',E'$ were found such that $E$ is isomorphic to $E'$ with $P'$ the image of $P$ under that isomorphism, then we list only one pair $P,E$. Further, for any $P,E$ pair corresponding to $h_2,h_3,h_4$ there is (by Proposition \ref{conjugatefactors}) another pair $P',E'$ corresponding to $\overline{h_2},\overline{h_3},\overline{h_4}$ (not necessarily in $S^+\times S \times S$) with $\hat{h}(P')=\hat{h}(P)$ but $E'$ not necessarily isomorphic to $E$; again, we list only one example, since $P',E'$ are then easily recovered from $h_2,h_3,h_4$ if desired. A full table has been made available online at \cite{Taylor} and will hopefully continue to grow: contributions are invited.

\begin{table}[h]
\begin{center}
{\small
\begin{tabular}{|p{5em}|c|p{26em}|l|}
\hline
$\hat{h}(P)$  & $\omega$ & $P$,\,$E$ & $u_2,u_3,u_4$\\
\hline
0.0038563\ldots & $\sqrt{3}$ &  $(497\omega+860,7523\omega+13029)$ on $y^2+\omega xy + y =x^3-\omega x^2 -(1544879\omega+2675806)x +1407381720\omega +2437656645$ & $1,\omega-1,2\omega-2$\\
\hline
0.0047223\ldots &$\frac{1}{2}+\frac{\sqrt{-7}}{2}$ & $(2\omega-4,3\omega+6)$ on $y^2+xy=x^3-2\omega x^2 -(\omega+2)x +5\omega+2$ & $\omega+1,-(2\omega+2),-(2\omega+10)$\\
\hline
0.0053416\ldots & $\sqrt{3}$ & $(1,0)$ on $y^2+xy+\omega y=x^3-\omega x^2 -5x+\omega+4$&$\omega+1,-(2\omega+2),-(8\omega+16)$\\
\hline
0.0054424\ldots &$\frac{1}{2}+\frac{\sqrt{-7}}{2}$ & $(\omega-2,5-4\omega)$ on $y^2+xy-\omega y=x^3+(3\omega-2)x^2 +(7\omega-5)x -12\omega-29$ & $1-\omega,6-2\omega,8\omega-24$\\
\hline
0.0058010\ldots &$\frac{1}{2}+\frac{\sqrt{-7}}{2}$ & $(2,-(3\omega+1))$ on $y^2+\omega xy+2\omega y=x^3-2x -\omega+3$ & $\omega,4-2\omega,16-8\omega$\\
\hline
0.0060112\ldots &$\frac{1}{2}+\frac{\sqrt{-7}}{2}$ & $(1,-2)$ on $y^2+\omega xy+y=x^3+x^2 -\omega x-\omega$ & $1,\omega-1,\omega+1$\\
\hline
0.0061272\ldots & $\sqrt{2}$ & $(1,0)$ on $y^2+xy+(1+\omega)y=x^3+\omega x^2-x-\omega$&$\omega+2,4\omega+6,20\omega+28$\\
\hline
0.0064724\ldots  & $\sqrt{-2}$ & $(2\omega,2\omega+1)$ on $y^2+\omega y =x^3-\omega x^2 +(2\omega+6)x + (\omega-3)$&$\omega,2\omega-4,-(4\omega+16)$\\
\hline
0.0069470\ldots  & $\sqrt{2}$ & $(\omega,-1)$ on $y^2+(\omega+1)xy+\omega y = x^3 -\omega x^2 - (\omega+3)x +1+\omega$&$1,\omega,\omega$\\
\hline
0.0072803\ldots &$\frac{1}{2}+\frac{\sqrt{-7}}{2}$ & $(0,2-\omega)$ on $y^2+xy=x^3-2\omega x^2 +(7\omega-6)x +2 -3\omega$ & $1-\omega,-(2\omega+2),24-8\omega$\\
\hline
0.0073349\ldots &$\frac{1}{2}+\frac{\sqrt{-7}}{2}$ & $(-(\omega+3),4\omega+2)$ on $y^2+xy+(\omega-1)y=x^3+(\omega+1)x^2 -(7\omega+1)x +\omega-11$ & $\omega,4-2\omega,8\omega$\\
\hline
0.0073479\ldots &$\frac{1}{2}+\frac{\sqrt{-7}}{2}$ & $(-\omega,\omega-1)$ on $y^2+(\omega-1)xy-\omega y=x^3+3\omega x^2 +(\omega-5)x +1-\omega$ & $\omega,\omega+2,5\omega+2$\\
\hline
0.0074870\ldots  & $\frac{1}{2}+\frac{\sqrt{-3}}{2}$ & $(1,0)$ on $y^2+(1-\omega)xy+ (1-\omega)y =x^3-2\omega x^2 +(2\omega-1)x$&$2-2\omega,-4,16$\\
\hline
0.0074943\ldots  & $\sqrt{-1}$ &  $(\omega,1)$ on $y^2+xy+\omega y=x^3+x^2+2x+\omega+2$&$\omega,\omega-1,-2$\\
\hline
0.0076951\ldots &$\frac{1}{2}+\frac{\sqrt{-7}}{2}$ & $(3-3\omega,2-2\omega)$ on $y^2+xy+(1-\omega)y=x^3+(3\omega-1)x^2 +(3\omega-1)x +3\omega-9$ & $1-\omega,-2,4$\\
\hline
0.0080799\ldots & $\frac{1}{2}+\frac{\sqrt{5}}{2}$ & $(1,2\omega-1)$ on $y^2+(1-\omega)xy+(1-\omega)y=x^3+\omega x^2+(\omega-2)x$&$2\omega,4\omega+4,32\omega+16$\\
\hline
0.0087764\ldots & $\sqrt{3}$ &$(0,0)$ on  $y^2+(1+\omega)y=x^3+(3+\omega)x^2+(2+2\omega)x$&$\omega+1,2\omega+2,4\omega+4$\\
\hline
0.0087786\ldots &$\sqrt{-1}$ & $(1,-\omega)$ on $y^2+\omega xy +y =x^3-(\omega+1)x^2$&$1-\omega,-2\omega,-4$\\
\hline
0.0088447\ldots & $\sqrt{2}$ & $(1-\omega,1)$ on $y^2+xy+(\omega+1)y=x^3+(\omega-1)x^2-(2\omega+2)x+1$&$\omega,-\omega,-2$\\
\hline
0.0089008\ldots &$\frac{1}{2}+\frac{\sqrt{-7}}{2}$ & $(2,-3)$ on $y^2+\omega xy=x^3-2\omega x^2 +\omega x +1$ & $2,2\omega-4,4\omega+8$\\
\hline
0.0089933\ldots & $\frac{1}{2}+\frac{\sqrt{5}}{2}$ & $(\omega+3,4\omega)$ on $y^2=x^3+\omega x^2 - (21\omega+16)x + 61\omega +41$&$\omega,3\omega+2,13\omega+8$\\
\hline
0.0089933\ldots & $\frac{1}{2}+\frac{\sqrt{5}}{2}$ & $(2\omega+3,3\omega+4)$ on $y^2+2\omega y=x^3-\omega x^2 - (6\omega+11)x + 22\omega +21$&$\omega+1,3\omega+2,21\omega+13$\\
\hline
0.0090543\ldots  & $\frac{1}{2}+\frac{\sqrt{-3}}{2}$ & $(2\omega-4,4\omega+4)$ on $y^2=x^3+(1-2\omega)x^2-12x+36\omega-12$&$\omega+1,-(3\omega+3),-(9\omega+9)$\\
\hline 
0.0091282\ldots & $\frac{1}{2}+\frac{\sqrt{5}}{2}$  & $(5,2-2\omega)$ on $y^2+(1-\omega)xy+(\omega-1)y = x^3 -x^2 - (34\omega+45)x + 158\omega + 149$&$1-\omega,-2\omega,-(4\omega+2)$\\ 
\hline
0.0091781\ldots & $\sqrt{3}$ & $(0,0)$ on $y^2+\omega y = x^3+\omega x^2$&$\omega,3\omega,-9\omega$\\
\hline
0.0093444\ldots & $\sqrt{2}$ & $(1,-1)$ on $y^2+(\omega+1)xy+y=x^3-x^2-(3\omega+4)x+2\omega+3$&$\omega,-2\omega,8\omega+8$\\
\hline
0.0097150\ldots & $\sqrt{2}$ & $(\omega,\omega)$ on $y^2+xy+\omega y=x^3-\omega x^2 +(\omega-2)x +2\omega +4$&$2,4\omega,-16$\\
\hline
0.0097217\ldots & $\frac{1}{2}+\frac{\sqrt{5}}{2}$  & $(\omega-1,4\omega)$ on $y^2=x^3+x^2-(15\omega+12)x + 27\omega+20$&$\omega,-(\omega+1),3\omega+2$\\
\hline
0.0097259\ldots & $\frac{1}{2}+\frac{\sqrt{5}}{2}$  & $(2\omega,4\omega)$ on $y^2=x^3+(1-\omega)x^2-(20\omega+16)x+76\omega+48$&$\omega,2\omega+1,-(13\omega+8)$\\
\hline
0.0097259\ldots & $\frac{1}{2}+\frac{\sqrt{5}}{2}$  & $(-\omega,4\omega+4)$ on $y^2=x^3+(1-\omega)x^2-(11\omega+13)x+27\omega+22$&$\omega+1,-(5\omega+3),-(21\omega+13)$\\
\hline
\end{tabular}
}
\end{center}
\caption{Some points $P$ on curves $E$ over quadratic fields $\Q(\sqrt{D})$ with height $\hat{h}(P)$ at most $0.01$.}
\label{lowheighttable}
\end{table} 

\end{document}